----------
X-Sun-Data-Type: default
X-Sun-Data-Name: torquot.tex
X-Sun-Charset: us-ascii
X-Sun-Content-Lines: 1024

\documentclass[12pt,a4paper,leqno]{article}

\usepackage{theorem}
\usepackage{epsfig}
\usepackage{amsfonts}

\theoremstyle{change}{\theorembodyfont{\slshape}
\newtheorem{theorem}{Theorem}[section]
\newtheorem{proposition}[theorem]{Proposition}
\newtheorem{lemma}[theorem]{Lemma}
\newtheorem{remark}[theorem]{Remark}
\newtheorem{corollary}[theorem]{Corollary}

}

{
\theorembodyfont{\rmfamily}
\newtheorem{example}[theorem]{Example}
\newtheorem{examples}[theorem]{Examples}
\newtheorem{definition}[theorem]{Definition}
}

\def\proof{\noindent{\bf Proof.}\enspace}
\def\endproof{ \quad $\kasten$}

\def\tq{\mathop{\lower 4pt \hbox{$ {\buildrel{/} \over {\scriptstyle{\rm tor}}} $}}\nolimits}
\def\cq{\mathop{\lower 4pt \hbox{$ {\buildrel{/} \over {\scriptstyle{\rm cat}}} $}}\nolimits}
\def\gq{\mathop{\lower 4pt \hbox{$ {\buildrel{/} \over {\scriptstyle{\rm good}}} $}}\nolimits}
\def\aq{\mathop{\lower 4pt \hbox{$ {\buildrel{/} \over {\scriptstyle{\rm aff}}} $}}\nolimits}

\def\mal{\! \cdot \!}

\def\h#1{\widehat{#1}}
\def\t#1{\widetilde{#1}}
\def\b#1{\overline{#1}}
\def\CC{{\mathbb C}}
\def\ZZ{{\mathbb Z}}
\def\RR{{\mathbb R}}

\def\mal{\mathbin{\! \cdot \!}}

\def\Lin{\mathop{\rm Lin}\nolimits}

\def\star{\mathop{\hbox{\rm star}}}
\def\conv{\mathop{\hbox{\rm conv}}}
\def\cone{\mathop{\hbox{\rm cone}}}

\def\Hom{\mathop{\rm Hom}\nolimits}

\def\osubset{\subset \kern-10pt {\rm o}\ } 

\def\topto#1{\mathop{\longrightarrow}\limits^{#1}}

\def\kasten{\mathord{\vbox{\hrule
                     \hbox{\vrule
                     \hskip5pt
                     \vrule height5pt
                     \vrule}
                     \hrule}}}
\def\text#1{\hbox{\rm #1}}

\def\bigtopmapright#1{\smash{\mathop{\hbox to 35pt{\rightarrowfill}}\limits^{#1}}}
\def\topmapright#1{\smash{\mathop{\hbox to 30pt{\rightarrowfill}}\limits^{#1}}}

\def\bigbotmapright#1{\smash{\mathop{\hbox to 35pt{\rightarrowfill}}\limits_{#1}}}
\def\botmapright#1{\smash{\mathop{\hbox to 30pt{\rightarrowfill}}\limits_{#1}}}

\def\bigtopmapleft#1{\smash{\mathop{\hbox to 35pt{\leftarrowfill}}\limits^{#1}}}
\def\topmapleft#1{\smash{\mathop{\hbox to 30pt{\leftarrowfill}}\limits^{#1}}}

\def\bigbotmapleft#1{\smash{\mathop{\hbox to 35pt{\leftarrowfill}}\limits_{#1}}}
\def\botmapleft#1{\smash{\mathop{\hbox to 30pt{\leftarrowfill}}\limits_{#1}}}

\def\rmapdown#1{\Big\downarrow\rlap{$\vcenter {\hbox{$\scriptstyle#1$}}$}}

\def\lmapdown#1{\llap{$\vcenter{\hbox {$\scriptstyle#1$}}$}\Big\downarrow}

\def\rmapne#1{\nearrow\rlap{\hbox{$\scriptstyle#1$}}}

\def\rmapsw#1{\swarrow\rlap{\hbox{$\scriptstyle#1$}}}

\def\lmapse#1{\llap{\hbox{$\scriptstyle#1$}}\searrow}

\def\ltexindent#1{\hbox to \hangindent{#1\hss}\ignorespaces}

\begin{document}

\thispagestyle{empty}

\begin{center}

{\LARGE\bf Quotients of Toric Varieties}

\medskip

{\LARGE\bf by the Action of a Subtorus}

\renewcommand{\thefootnote}{\fnsymbol{footnote}}

\bigskip

Annette A'Campo--Neuen %
and %
J\"urgen Hausen%

\bigskip

\end{center}

\begin{abstract}\noindent
We consider the action of a subtorus  of the big torus  on a toric variety. The aim of the paper is to define a natural notion of a quotient for this setting and to give an explicit algorithm for the construction of this quotient from the combinatorial data corresponding to the pair consisting of the subtorus and the toric variety. Moreover, we study the relations of such quotients with good quotients. We construct a good model, i.e. a dominant  toric morphism from the given toric variety to some ``maximal'' toric variety having a good quotient by the induced action of the given subtorus.
\end{abstract}

\section*{Introduction}

Let $X$ be an algebraic variety with a regular action of an algebraic
group $G$. A categorical quotient is a morphism $p: X \to Y$ which is
$G$--invariant (i.e.~constant on $G$--orbits) and satisfies the
following universal property: every other $G$--invariant morphism $f:
X \to Z$ factors uniquely through $p$ (see \cite{Mu}).

Though this universal property seems to be a minimal requirement for a
quotient, there is no hope for the general existence of categorical
quotients. (See e.g. \cite{AC;Ha} for an explicit example of a
$\CC^*$--action on a smooth four-dimensional toric variety which does
not have a categorical quotient, even if one allows the quotient space
$Y$ to be an algebraic or analytic space.)

In the present article we consider toric varieties $X$ with an action
of an algebraic torus $H$; we refer to these varieties as toric
$H$-varieties. The specialization of the definition of the categorical
quotient to the category of toric varieties leads to the following
notion: We call a  toric morphism $p: X \to Y$ a {\it toric quotient},
if it is $H$-invariant and every $H$-invariant toric morphism  factors
uniquely through $p$. For this kind of quotient we can actually prove
the existence (see Theorem \ref{torquot}):

\medskip

{\it For every toric $H$-variety $X$ there exists a toric quotient.}

\medskip

Our proof of this result is constructive. In fact, we introduce the
notion of a {\it quotient fan} of a fan by some sublattice (see
Section 2) and give an algorithm for the calculation of this quotient
fan. We obtain the existence of toric quotients by applying this
algorithm to the fan $\Delta$ of $X$ and the lattice $L$ of one
parameter subgroups of the acting torus $T$ of $X$ factoring through
$H$.

A particularly important notion of quotient is the so-called good quotient (see \cite{Se}) generalizing the quotients occuring in Mumford's geometric invariant theory for projective varieties. Unfortunately, good quotients exist only under very special circumstances. However, for any toric $H$-variety $X$  we can construct a {\it good model} $\b{X}$. More precisely, we show (see Theorem \ref{goodmod}):

\medskip

{\it There exists a ``maximal'' toric $H$-variety $\b{X}$ with a good quotient such that there is a dominant $H$-equivariant toric morphism from $X$ to $\b{X}$. The good quotient of $\b{X}$ by $H$ coincides with the toric quotient of $X$ by $H$.}

\medskip

In fact, the good model  defines an adjoint functor to the forgetful functor from the category of toric $H$-varieties with good quotients into the category of toric $H$-varieties. Again our proof of the existence of the good model is constructive and works in terms of fans. The good model can be used to understand the obstructions for the existence of a good quotient.

\medskip

The authors would like to thank G.~Barthel, A.~Bia\l ynicki-Birula,  L.~Kaup and J.~\'Swi\c{e}cicka for their interest in the subject and for many helpful discussions.

\section{Toric Quotients}

First we briefly recall some of the basic definitions. A normal
algebraic variety $X$ is called a {\it toric variety} if there is an
algebraic action of a torus $T$ on $X$ with an  open orbit. We always
assume the action to be effective and refer to $T$ as the {\it acting
  torus} of $X$. For every toric variety $X$ we fix a point  $x_0$ in
its open orbit which we call the base point of $X$.

Let $X$, $X'$ be toric varieties with acting tori $T$, $T'$ and base
points $x_0$ and $x_0'$, respectively. A regular map $f : X \to X'$ is
called a {\it toric morphism} if $f(x_0)=x_0'$ and there is a
homomorphism $\varphi: T \to T'$ such that $f(t \mal x) = \varphi(t)
\mal f(x)$ for every $(t,x) \in T \times X$.

Now let $H$ be any algebraic torus. We call a given toric variety $X$
with acting torus $T$ a {\it toric $H$-variety}, if $H$ acts on $X$ by
means of a morphism $H \times X \to X$, $(h,x) \mapsto h \! * \! x$ of
algebraic varieties such that the actions of $H$ and $T$ on $X$
commute.

\begin{remark}\label{actionhom}
If $X$ is a toric $H$-variety, then there is a homomorphism $\psi$
from $H$ into the acting torus $T$ of $X$ such that the action of $H$
on $X$ is given by $h \! * \! x = \psi(h) \mal x$.
\end{remark}

\proof The action of $H$ permutes the $T$-orbits since it commutes with the 
$T$-action.  The open orbit $T\mal x_0$ is even $H$-stable because there is 
only one $T$-orbit of maximal dimension. Since the action of $T$ is effective,
for every $h\in H$, there is a unique element $\psi(h)$ in $T$ such that 
$h \! * \! x_0=\psi(h)\mal x_0$. Now it is straightforward to check that the
map $H\to T$, $h\mapsto \psi(h)$ has the required properties.
\endproof

\bigskip

An $H$-equivariant  toric morphism  $f : X \to X'$ of two toric $H$-varieties 
will  be called a {\it toric $H$-morphism}. If the action of $H$ on $X'$ is 
trivial, which means that $f$ is constant on $H$-orbits, then we will
say that $f$ is {\it $H$-invariant}.

\begin{definition} We call an $H$-invariant toric morphism $p : X \to
  Y$ a {\it toric quotient} for the toric $H$-variety $X$, if it has
  the following universal property: for every $H$-invariant toric
  morphism $f : X \to Z$ there is a unique toric morphism $\tilde f: Y
  \to Z$ such that the diagram

\smallskip

$$\matrix{%
X \qquad \topto{f} \qquad Z \cr
\lmapse{p} \qquad \rmapne{\tilde f} \cr
   Y   \cr}%
$$

\smallskip

\noindent is commutative. Note that the toric quotient $Y$ is uniquely
determined by this property. We denote the quotient space $Y$ also by $X \tq H$.
\end{definition}

\begin{example}\label{C2} Let $H := \CC^*$ act on the toric variety $X
  := \CC^2$ by the homomorphism $t \mapsto (t^a, t^b)$, where $a$ and
  $b$ are relatively prime integers and $a > 0$. Then one can verify
  directly that the toric quotient of $X$ is the constant map $\CC^2
  \to \{0\}$ if $b > 0$, and that otherwise it is the morphism

\smallskip

$$ p : \CC^2 \to \CC, \quad (z,w) \mapsto \cases{ w & if $b=0$, \cr
  z^{-b} w^a & if $b < 0$. \cr}$$

\smallskip

So $p$ corresponds to the inclusion $\CC[X]^H \subset \CC[X]$,
i.e. $p$ equals the categorical quotient for the action of $H$ on the
affine variety $X$. In fact, this holds generally for affine toric
$H$-varieties (see Example \ref{affquot}). \end{example}

Further basic examples of toric quotients are invariant toric
fibrations, e.g. line bundles on toric varieties. For the toric
quotient defined here we have not only uniqueness but  also the
existence.

\begin{theorem}\label{torquot}
Every toric $H$-variety $X$ has the toric quotient $p : X \to X \tq
H$.
\end{theorem}

For the proof of this theorem we use the description of toric
varieties by means of  fans. Let us first fix some notation. For an
algebraic torus $T$, denote by $N_T$ the lattice $\Hom(\CC^*, T)$ of
its one parameter subgroups. A fan $\Delta$ in $N_T$ is a finite set
of strictly convex rational polyhedral cones in $N_T^\RR := N_T
\otimes_{\ZZ} \RR$ satisfying the following two conditions: any two
cones of $\Delta$ intersect in a common face, and if $\sigma \in
\Delta$, then $\Delta$ also contains all the faces of $\sigma$. We
denote a fan $\Delta$ in $N_T$ also as a pair $(N_T, \Delta)$.

For every fan $\Delta$ in $N_T$, there is a corresponding toric
variety $X_\Delta$ with the acting torus $T$ (as basic references for
this construction, see e.g. \cite{Fu} and \cite{Od}). The assignment
$\Delta \mapsto X_\Delta$ yields an equivalence between the category of
fans and the category of toric varieties (with fixed base point),
where maps of fans correspond to toric morphisms. Recall that a map of
fans $F: (N,\Delta) \to (N',\Delta')$ by definition is a $\ZZ$-linear
homomorphism from $N$ to $N'$, also denoted by $F$, such that for
every cone $\sigma\in\Delta$ there is a cone $\tau\in\Delta'$ with
$F^{\RR}(\sigma)\subset \tau$ ($F^{\RR}:N^{\RR}  \to {N'}^{\RR}$ is
the scalar extension of $F$).

Now, if a torus $H$ acts on a toric variety $X_\Delta$ by a
homomorphism $\varphi$ from $H$ to the acting torus $T$ of $X_\Delta$,
let $L$ denote the (primitive) sublattice of $N_T$ corresponding to
the subtorus $\varphi(H)$ of $T$. Then a  toric morphism $f : X_\Delta
\to X_{\Delta'}$ is $H$-invariant if and only if the corresponding map
of fans $F : (N,\Delta) \to (N',\Delta')$ satisfies $L \subset
\ker(F)$. So in the language of fans, Theorem \ref{torquot} reads as
follows:

\begin{theorem}\label{torquot1} Let $\Delta$ be a fan in a lattice $N$
  and let $L$ be a primitive sublattice of $N$. Then there is a map of
  fans $P:(N,\Delta)\to (\t{N},\t{\Delta})$ with $L \subset \ker(P)$
  such that the following universal property is satisfied: for every
  map of fans $F : (N,\Delta) \to (N',\Delta')$ with $L\subset \ker F$
  there is a unique map of fans $\t{F} : (\t{N},\t{\Delta}) \to
  (N',\Delta')$ with $F = \t{F} \circ P$.
\end{theorem}

The fan $\t{\Delta}$ occuring in the above theorem will be called the
{\it quotient fan} of $\Delta$ by $L$. Note that our concept of a
quotient fan differs from the notion introduced in \cite{KaStZe},
since we require the existence of a map of fans  from $\Delta$ to
$\t{\Delta}$. We will prove Theorem \ref{torquot1} in the next section
by describing an explicit algorithm to construct the quotient fan. The
algorithm starts with projecting cones of $\Delta$ to $N/L$. But then
two types of difficulties  occur:

Firstly, the projected cones in general are no longer strictly
convex. Secondly, it can happen that the projected cones do not
intersect in a common face. Therefore the construction requires an
iteration of steps refining the first na\"\i ve approach. The first of
the above-mentioned difficulties already occurs in Example
\ref{C2}. Here is its fan-theoretic version:

\begin{example}\label{C2fan}\begingroup\rm
The fan $\Delta$ of the toric variety $\CC^2$ consists of the faces of
the cone $\sigma \in \RR^2$ spanned by the canonical basis vectors
$e_1$ and $e_2$. The action of $H = \CC^*$ on $X_\Delta$ considered in
Example \ref{C2} corresponds to the line $L$ through the point
$(a,b)$.

Let $P : \ZZ^2 \to \ZZ^2 / L$ denote the projection. If $b \le 0$, the
quotient fan $\t{\Delta}$ of $\Delta$ by $L$ is the fan of faces of
$P^\RR(\sigma)$ in $\t{N} := \ZZ^2 / L$. If $b > 0$, then
$P^\RR(\sigma)$  fails to be strictly convex and the quotient fan is
just the zero fan in $\t{N} = \{0\}$.
\endgroup\end{example}

\section{Computation of the Quotient Fan}

Let $N$ be a lattice, i.e.~a free $\ZZ$-module of finite rank. In this
section we construct the quotient fan of a fan $\Delta$ in $N$ by a
primitive sublattice $L$ of $N$ and thereby prove Theorem
\ref{torquot1}. In fact our construction is done in a more general
framework. We will not only consider fans but also sets of convex
rational polyhedral cones which are not required to be strictly convex
nor to intersect pairwise in a common face.

More precisely, we will speak of a {\it system $S$ of $N$-cones} if
$S$ is a finite set of convex cones in the space $N^\RR :=
N\otimes_{\ZZ}\RR$ such that every $\sigma \in S$ is generated by
finitely many vectors of $N$. A {\it map} $F : (N,S) \to (N',S')$ of a
system $S$ of $N$-cones to a system $S'$ of $N'$-cones is a lattice
homomorphism from $N$ to $N'$, also denoted by $F$, such that for
every $\sigma \in S$ there is a cone $\tau\in S'$ with $F^{\RR}
(\sigma) \subset \tau$. This notion generalizes the concept of a map
of fans.

We also need the following ``intermediate'' notion: A system $\Sigma$
of $N$-cones is called a {\it quasifan} in $N$, if for each $\sigma \in
\Sigma$ the faces of $\sigma$ also belong to $\Sigma$ and for any two
cones $\sigma$ and $\sigma'$ of $\Sigma$ the intersection $\sigma \cap
\sigma'$ is a face of $\sigma$. So a quasifan is a fan if all its
cones are strictly convex. A map of two quasifans is just a map of the
underlying systems of cones.

\begin{definition}\label{defquotfan}
Let $N$ be a lattice and let $S$ be a system of $N$-cones. If $L
\subset \h{L} \subset N$ are primitive sublattices, then we call a
(quasi-) fan $\t{\Delta}$ in $\t{N} := N / \h{L}$ a {\it quotient
  (quasi-) fan} of $\Delta$ by $L$ if it has the following properties:

\begin{enumerate}
\item The projection $P : N \to \t{N}$ defines a map of the systems
  $S$ and $\t{\Delta}$ of cones.
\item For every map $F: (N,S) \to (N',\Delta')$ from $S$ to a (quasi-)
  fan $\Delta'$ in a lattice $N'$ with $F(L)= 0$, there is a map
  $\t{F}: (\t{N},\t{\Delta}) \to (N',\Delta')$ of (quasi-) fans such
  that the following diagram is commutative:

\smallskip

$$\matrix{%
(N,S) \quad \topto{F} \quad (N',\Delta') &\cr
\lmapse{P} \qquad \quad \rmapne{\t{F}} &\cr
   (\t{N},\t{\Delta})   &. \cr}%
$$

\end{enumerate}
\end{definition}

By definition, quotient fans and quotient quasifans are uniquely
determined. These two notions are related to each other by the
following:

\begin{remark}\label{quasifan} Let $\Sigma$ be a quasifan with maximal
  cones $\sigma_1, \ldots, \sigma_r$ in a lattice $N$. For the maximal
  sublattice $L$ of $N$ contained in $\bigcap_{i=1}^r \sigma_i$ let
  $P: N \to \t{N} := N/L$ denote the projection. Then the cones
  $P^\RR(\sigma_1), \ldots, P^\RR(\sigma_r)$ are the maximal cones of
  the quotient fan $\t{\Delta}$ of $\Sigma$ by $L$.
\end{remark}

\proof Set $\sigma_0 := \cap_{i=1}^r \sigma_i$. Then $\sigma_0$ is a
cone with $V := L^\RR$ as the smallest face. Since $\sigma_0$ is a
face of each $\sigma_i$, it follows that $V = \ker(P^\RR)$ is also the
smallest face of every $\sigma_i$. This implies
$(P^\RR)^{-1}(P^\RR(\sigma_i)) = \sigma_i$ for every $i$.

As a consequence we obtain that every cone $P^\RR(\sigma_i)$ is
strictly convex. Now we check that for any two $i$ and $j$ the cones
$P^\RR(\sigma_i)$ and $P^\RR(\sigma_j)$ intersect in a common
face. Note that

$$ P^\RR(\sigma_i) \cap P^\RR(\sigma_j) = P^\RR(\sigma_i \cap
\sigma_j).$$

Choose a supporting hyperplane $W$ of $\sigma_i$ defining the face
$\sigma_i\cap\sigma_j$. Since $W$ contains $V$, its projection
$P^\RR(W)=W/V$ is a supporting hyperplane of $P^\RR(\sigma_i)$ that
cuts out $P^\RR(\sigma_i \cap \sigma_j)$. Therefore $P^\RR(\sigma_i)
\cap P^\RR(\sigma_j)$ is a face of  $P^\RR(\sigma_i)$.

So the cones $P^\RR(\sigma_1), \ldots, P^\RR(\sigma_r)$ together with
their faces define a fan $\t{\Delta}$ in $\t{N}$. By construction,
$\t{\Delta}$ satisfies the properties of a quotient fan of $\Sigma$ by
$L$. \endproof

\bigskip

The main result of this section is the following:

\begin{theorem}\label{quotcalc}
For a given system $S$ of $N$-cones and a primitive sublattice $L$ of
$N$, there is an algorithm to construct the quotient fan $\t{\Delta}$
of $S$ by $L$.
\end{theorem}

\proof Set $N_1 := N/L$ and let $P_1 : N \to N_1$ denote the
projection. We first construct a quotient quasifan $\Sigma$ in $N_1$
of $S$ by $L$ by means of the following procedure:

\bigskip

{\it Initialization.}\enspace Define $S_1$ to be the system of
$N_1$-cones consisting of those of the $P_1^\RR(\sigma)$, $\sigma \in
\Delta$, that are maximal with respect to set-theoretic inclusion.

\bigskip

{\it Loop.}\enspace While there exist cones $\tau_1$ and $\tau_2$ in
$S_1$ such that $\tau_1 \cap \tau_2$ is not a face of $\tau_1$ do the
following: Let $\varrho_2$ be the minimal face of $\tau_2$ that
contains $\tau_1 \cap \tau_2$. If $\varrho_2 \not\subset \tau_1$,
replace $\tau_1$ by the convex hull $\conv(\tau_1 \cup \varrho_2)$ of
$\tau_1\cup \varrho_2$. Otherwise let $\varrho_1$ be the minimal face
of $\tau_1$ that contains $\tau_1 \cap \tau_2$ and replace $\tau_2$ by
$\conv(\tau_2 \cup \varrho_1)$. Omit all cones of $S_1$ that are
properly contained in the new one.

\bigskip

{\it Output.}\enspace Let $\Sigma$ be the system of $N_1$-cones
consisting of all the faces of the cones of $S_1$.

\bigskip

The above loop is finite: passing through the loop does not increase
the number $\vert S_1 \vert$ of cones of $S_1$. So, after finitely
many, say $K$, steps $\vert S_1 \vert$ stays fixed. For each iteration,
there is a cone $\tau$ of $S_1$ that is replaced by a strictly larger
cone of the form $\conv(\tau \cup \varrho)$ with a face $\varrho$ of
some other cone of $S_1$. According to Lemma \ref{gencone} below we
obtain
$$ \vert \tau \cap P_1^\RR(S^{(1)}) \vert \; < \; \vert \conv(\tau
\cup \varrho) \cap P_1^\RR(S^{(1)}) \vert, $$
where $S^{(1)}$ denotes a minimal set of generators of the cones of
$S$. Thus in every step after the first $K$ steps the number $
\sum_{\tau \in S_1} \vert \tau \cap P_1^\RR(S^{(1)}) \vert$
strictly increases. But this can only happen a finite number of
times. So the loop is indeed finite.

Now by construction $\Sigma$ is a quasifan. We have to verify that it
fullfills the property ii) of Definition \ref{defquotfan}. So let $F:
(N,S) \to (N', \Sigma')$ be a map of quasifans with $L \subset
\ker(F)$. Then there is a lattice homomorphism $F_1 : N_1 \to N'$ with
$F = F_1 \circ P_1$. Clearly $F_1$ defines a map from the system $S_1$
of cones defined in the initialization to the system $\Sigma'$ of
cones.

Assume that after $n$ iterations of the loop, $F_1$ still defines a
map of the systems of cones $S_1$ and $\Sigma'$, and that in the next
step we replaced the cone $\tau_1$ by $\conv(\tau_1 \cup \varrho_2)$,
where $\varrho_2$ is a minimal face of $\tau_2$ such that $\tau_1 \cap
\tau_2 \subset \varrho$. We have to check that there is a cone in
$\Sigma'$ containing $F_1^\RR(\conv(\tau_1 \cup \varrho_2))$. Let
$\tau_1'$ and $\tau_2'$ be cones of $\Sigma'$ such that
$F_1^\RR(\tau_1) \subset \tau_1'$ and $F_1^\RR(\tau_2) \subset
\tau_2'$. Then
$$ F_1^\RR(\varrho_2)^\circ \; \cap \; (\tau_1' \cap \tau_2') \; \ne
\; \emptyset. $$
Since $\tau_1' \cap \tau_2'$ is a face of $\tau_2'$ and
$F_1^\RR(\varrho_2) \subset \tau_2'$, we obtain $F_1^\RR(\varrho_2)
\subset \tau_1' \cap \tau_2'$. This implies $F_1^\RR(\varrho_2)
\subset \tau_1'$. In particular, it follows that
$$  F_1^\RR(\conv(\tau_1 \cup \varrho_2)) \subset \tau_1'.$$
Thus after $\tau_1$ is replaced by $\conv(\tau_1 \cup \varrho_2)$ the map
$F_1$ still defines a map of the systems $S_1$ and $\Sigma'$ of cones.

Repeating this argument we obtain that $F_1$ defines also a map of the
quasifans $\Sigma$ and $\Sigma'$. Thus $\Sigma$ fullfills the desired
universal mapping property and hence it is the quotient quasifan of
$S$ by $L$.

Now let $V$ denote the maximal linear subspace contained in the
intersection of all maximal cones of the quasifan $\Sigma$. Set $L_1
:= N_1 \cap V$. Then, according to Lemma \ref{quasifan}, the quotient
fan $\t{\Delta}$ in $N_1 / L_1$ of $\Sigma$ by $L_1$ is obtained by
projecting the maximal cones of $\Sigma$ to $N_1^\RR / L_1^\RR$. It
follows that with $\h{L} := P_1^{-1}(L_1)$, the fan $\t{\Delta}$ is
also the quotient fan of $S$ by $L$. \endproof

\bigskip

We have used the following elementary fact about cones:

\begin{lemma}\label{gencone} Let $\sigma=\cone(v_1,\dots,v_r)$ be the
  (not necessarily strictly) convex cone spanned by $v_1,\dots,v_r$ in
  a real vector space $V$. Then every face $\tau$ of $\sigma$ is
  generated as a cone by the vectors in $\tau \cap
  \{v_1,\dots,v_r\}$. \endproof
\end{lemma}

As a consequence of the construction of the quotient fan we note:

\begin{remark}\label{maxcones}
Let $\tau_1, \ldots, \tau_k$ be the maximal cones of the quotient fan of $S$ by $L$, and let $F(S)$ denote the set of all faces of the cones of $S$. Then
$$ \tau_i = \conv \left( \bigcup_{\sigma \in F(S); \; P^\RR(\sigma)
    \subset \tau_i} P^\RR(\sigma) \right). $$
\end{remark}

\noindent {\bf Proof of Theorems \ref{torquot} and \ref{torquot1}.}\enspace It
suffices to verify Theorem \ref{torquot1}. So, let $\t{\Delta}$ be the
quotient fan of $\Delta$ by $L$. Then, if $P$ denotes the projection
from $N$ onto the lattice of $\t{\Delta}$, we have only to check that
the factorization of every $L$-invariant  map of fans $F : (N,\Delta)
\to (N',\Delta')$ through $P$ is unique. But this follows from the
fact that $P$ is surjective by construction. \endproof

\bigskip

In the case of small codimension of $L$ in $N$, there is an easy
explicit description of the quotient fan:

\begin{example}Let $\Delta$ be a fan in a lattice $N$, and let $L
  \subset N$ be a primitive sublattice of codimension $2$. Denote by
  $P$ the canonical projection $N \to N/L$ and define an equivalence
  relation on the set of maximal cones of $\Delta$ as follows:

Set $\sigma \sim \tau$ if there is a sequence $\sigma=\sigma_0$,
$\sigma_1, \ldots, \sigma_r = \tau$ of cones $\sigma_i \in \Delta$
such that $P^\RR(\sigma_i^\circ) \cap P^\RR(\sigma_{i+1}^\circ) \not=
\emptyset$. For each maximal cone $\sigma\in\Delta$ denote by
$\b{\sigma}$ the convex hull of the union of all maximal cones
$\tau\sim\sigma$.

Let $V$ denote the sum of all linear subspaces of the cones
$P^\RR(\b{\sigma})$, $\sigma \in \Delta$ and set $\h{L} := N \cap
P^{-1}(V)$. Moreover, let $P : N \to N / \h{L}$ be the
projection. Then the faces of the cones $Q^\RR(\b{\sigma})$, where
$\sigma$ varies over the maximal cones of $\Delta$, form the quotient
fan of $\Delta$ by $L$.
\end{example}

In \cite{Ew} a special case of our notion of the quotient fan is
introduced for the abstract description of orbit closures of the
acting torus of a toric variety. In fact these orbit closures are
toric quotients of certain neighborhoods:

\begin{example}\label{Bahnabschluesse}
Let $\Delta$ be a fan in a lattice $N$. For a cone $\tau \in \Delta$,
let $x_\tau$ be the corresponding distinguished point in the toric
variety $X_\Delta$ (see [Fu, p. 27]). Let $B_\tau$ be the orbit of the
acting torus $T$ of $X_\Delta$ through $x_\tau$. Denoting by
$\star(\tau)$ the set of all cones $\sigma \in \Delta$ that contain
$\tau$ as a face, we obtain the closure of the orbit $B_\tau$ as

$$ V(\tau) := \overline{B_\tau} = \bigcup_{\sigma \in \star(\tau)}
B_{\sigma}.$$

The union $U(\tau)$ of the affine charts $X_{\sigma}$, $\sigma \in \star(\tau)$, is an open $T$-invariant neighbourhood of the orbit closure $V(\tau)$. For the set of maximal cones of the fan $\Delta(\tau)$ corresponding to $U(\tau)$ we have
$$\Delta(\tau)^{\max} = \Delta^{\max} \cap \star(\tau).$$

Let $L$ be the intersection of the linear hull $\Lin(\tau)$ 
of $\tau$ in $N^{\RR}$ with the lattice $N$, and let $P : N \to N/L$ denote the
projection. Then the cones $P^{\RR}(\sigma)$, $\sigma \in
\Delta(\tau)^{\max}$ are the maximal cones of the quotient fan
$\t{\Delta}(\tau)$ of $\Delta(\tau)$ by $L$. Moreover
$\t{\Delta}(\tau)$ is the fan of $V(\tau)$, viewed as a toric variety
with acting torus $T/T_{x_\tau}$ (see e.g. [Fu, p. 52]). In other
words, the toric morphism $p : U(\tau) \to V(\tau)$ associated to $P$
is the toric quotient of $U(\tau)$ by $T_{x_\tau}$.
\end{example}

\section{Good Models}

Let $X$ be an algebraic variety with a regular action of a reductive
group $G$. If $X$ is affine, then the categorical quotient for this
action always exists, and is given by the morphism corresponding to
the inclusion of the algebra $\CC[X]^G$ of $G$-invariant regular
functions on $X$ into $\CC[X]$. For general $X$, the idea of glueing
affine quotients of $G$-stable affine charts leads to the following
definition (see \cite{Se}):

A $G$-invariant morphism $p : X \to Y$ of algebraic varieties is
called a {\it good quotient}, if there exists a covering $(U_i)_{i \in
  I}$ of $Y$ by affine open sets such that every $W_i := p^{-1}(U_i)$
is affine and the restriction $p_{\vert W_i} : W_i \to U_i$ is the
categorical quotient for the action of $G$ restricted to $W_i$. If in
addition the morphism $p$ separates orbits, it is called a {\it
  geometric quotient}.

Now, coming back to the setting of toric $H$-varieties, we will first
give the description of the affine case in terms of fans:

\begin{example}\label{affquot}Let $T$ be an algebraic
  torus and let $\sigma$ be a rational strictly convex cone in
  $N_T^\RR$. Denote by $X_\sigma$ the associated affine toric
  variety. For a given subtorus $H \subset T$ let $L$ be the
  sublattice of $N_T$ corresponding to $H$. Define $\tau$ to be the
  maximal face of $\sigma$ with $L \cap \tau^\circ \ne \emptyset$ and
  set
$$ \h{L} :=  (L^{\RR} + \Lin(\tau)) \cap N_T.$$
Denote by $P : N_T \to N_T / \h{L}$ the canonical projection. Then
$P^\RR(\sigma)$ is a rational strictly convex cone in $N_T^{\RR} /
\h{L}^{\RR}$, and  the toric morphism $p : X_\sigma \to
X_{P^\RR(\sigma)}$ associated to $P$ is the toric quotient for the
action of $H$ on $X$.

The coordinate algebra of $X_{P^\RR(\sigma)}$ can be identified with
the algebra $\CC[X]^H$ of $H$-invariant regular functions on $X$,
since every $H$-invariant character of $T$ extending to a regular
function on $X$ factors through $p$. This shows that $p$ is also the
categorical quotient.
\end{example}

If a toric $H$-variety $X$ has a good quotient $p: X\to Y$, then it
follows that $Y$ is a toric variety and $p$ is a toric
morphism. Moreover, we can conclude that if a good quotient exists, it
coincides with the toric quotient. Conversely, as a consequence of
Example 3.1, our procedure for the calculation of the quotient fan
yields a good quotient if and only if it produces an affine map. So
we can characterize fan theoretically when a given toric quotient is
good (see also [Sw]):

\begin{proposition}\label{goodcrit}
Suppose $p : X_\Delta \to X_{\t{\Delta}}$ is the toric quotient of a toric 
$H$-variety $X_\Delta$. Let $P : (N_T,\Delta) \to (N_{\t{T}},\t{\Delta})$ be the associated map of fans. Then $p$ is good if and only if the following two conditions are satisfied:
\begin{enumerate}
\item For every maximal cone $\tau_i \in \t{\Delta}$ there is a
  maximal cone $\sigma_i \in \Delta$ such that $P^{\RR}(\sigma_i) =
  \tau_i$.
\item Every ray $\varrho \in \Delta^{(1)}$ with
  $P^\RR(\varrho)\subset\tau_i$ is contained in $\sigma_i$.
\end{enumerate}
Moreover, $p$ is geometric, if in addition $\dim\tau_i=\dim\sigma_i$
for all $i$. \endproof
\end{proposition}

Good quotients have excellent properties, but unfortunately they only
rarely exist. Bia\l ynicki-Birula and \'Swi\c
ecicka (\cite{joanna}) give a complete description of all open  subsets of $X$ having
a good quotient. Instead of looking at subsets one can also try to
modify $X$ to obtain a toric $H$-variety having a good quotient.  This
approach leads to the following notion:

\begin{definition} Let $p : X \to X \tq H$ denote the toric quotient
  of the action of $H$ on $X$. Suppose that $g: X \to \b{X}$ is a
  dominant toric $H$-morphism to a toric $H$-variety $\b{X}$ having a
  good quotient. Then we call  $g$ a {\it good model} for the toric
  $H$-variety $X$, if it has the following universal property: If $f :
  X \to Z$ is a toric $H$-morphism and the toric $H$-variety $Z$ has a
  good quotient, then there is a unique toric $H$-morphism $\b{f} :
  \b{X} \to Z$ such that the following diagram is commutative:

\smallskip

$$
\matrix{%
X \qquad \topto{g} \qquad \b{X} &\cr
\lmapse{f} \qquad \rmapsw{\b{f}} &\cr
   Z   &. \cr}
$$
\end{definition}

\medskip

Being defined by a universal property, a good model is unique up to
isomorphism. If $g$ is a good model, then there is a unique toric
morphism $\b{p}:\b{X}\to X\tq H$ such that the diagram

\smallskip

$$\matrix{%
X \qquad \topto{g} \qquad \b{X} &\cr
\lmapse{p} \qquad \rmapsw{\b{p}} &\cr
X\tq H       &\cr}%
$$ 

\medskip

\noindent is commutative. It follows that $\b{p}$ is in fact the toric
and hence the good quotient for the action of $H$ on $\b{X}$. Before
proceeding to the general construction of the good model we  give some
elementary examples.

\begin{examples}\label{C2var}
\enspace{\bf a)\enspace} Let $X$ be  $\CC^2\backslash \{0\}$ and let
$H$ be the subtorus $\{(t,t^{-1});t\in\CC^*\}$ of the acting torus
$(\CC^*)^2$ of $X$. Then the toric quotient is the map $p:X\to \CC$
defined by $p(z,w) = zw$ and the good model is just the inclusion of
$X$ in $\CC^2$ (compare Example \ref{C2}). So in this case the
``missing'' fixed point $0$ has to be added to $X$.

\smallskip

{\bf b)\enspace} Let $X$ be the blow-up of $\CC^2$ at the point
$0$. Then the action of the  torus $H$ in a) as well as the toric
quotient map extend naturally to $X$. The good model of $X$ is the
blow-down map $g: X \to \CC^2$ contracting the exceptional curve to a
point. 

\smallskip

{\bf c)\enspace} If $X$ is complete, then the toric quotient space is
also complete and the good model equals the toric quotient.
\end{examples}

The main result of this section is the following:

\begin{theorem}\label{goodmod}
Every toric $H$-variety $X$ has a good model. If $X = X_\Delta$, then the good model is obtained as follows: Let $P : (N,\Delta) \to (\t{N},{\t{\Delta}})$ be the map of fans corresponding to the toric quotient $p : X_\Delta \to X_{\t{\Delta}}$ of the action of $H$ on $X_\Delta$. For every maximal cone $\tau_i$, $i = 1, \ldots, r$, of $\t{\Delta}$, set
$$ {\sigma}_i := \conv \{\varrho \in \Delta^{(1)}; \; 
P^\RR(\varrho) \subset \tau_i\}  .$$
Moreover, let $V$ be the maximal linear subspace contained in $\bigcap_{i=1}^r \sigma_i$, set $L := V \cap N$ and let $G : N \to \b{N} := N/L$ denote the projection. Then  $G(\sigma_1),\dots,G(\sigma_r)$  are the maximal cones of a fan $\b{\Delta}$ in $\b{N}$, the projection $G$ defines a map of fans from $\Delta$ to $\b{\Delta}$ and the associated toric morphism $g : X_{\Delta} \to X_{\b{\Delta}}$ is the good model for $X_\Delta$.
\end{theorem}

The assignment $X \mapsto \b{X}$ is even functorial. More precisely, if $X$ and $X'$ are toric $H$-varieties with good models $g : X \to \b{X}$ and $g' : X' \to \b{X'}$, then for every toric $H$-morphism $f : X \to X'$, there is a unique toric $H$-morphism $\b{f} : \b{X} \to \b{X'}$ such that $\b{f} \circ g = g' \circ f$. A fancy formulation of the properties of the good model in the language of categories is the following:

\begin{corollary}
The assignment $X \mapsto \b{X}$ is adjoint to the forgetful functor from the category of toric $H$-varieties with good quotients into the category of toric $H$-varieties.
\end{corollary}

\bigskip

\noindent{\bf Proof of Theorem \ref{goodmod}.}\enspace First we prove that $\sigma_1,\dots,\sigma_r$ are the maximal cones of a quasifan $S$. Let ${\sigma}_i$ and ${\sigma}_j$ be two cones of $S$ and let $\sigma$ denote the minimal face of $\sigma_i$ containing the intersection $\sigma_i \cap \sigma_j$. Then there is a vector $v \in (\sigma_i \cap \sigma_j)^\circ \cap \sigma^\circ$. Moreover, for this $v$ we have
$$ P^\RR(v) \; \in \; P^\RR((\sigma_i\cap \sigma_j)^\circ \cap \sigma^\circ) \; \subset \; P^\RR((\sigma_i\cap \sigma_j)^\circ) \cap P^\RR(\sigma^\circ)
 \; \subset \; (\tau_i\cap\tau_j) \cap P^\RR(\sigma)^\circ.$$
In particular, the intersection of $P^\RR(\sigma)^\circ$ with $\tau_i\cap\tau_j$ is not empty. Since $\tau_i\cap\tau_j$ is a face of $\tau_i$ and the cone $P^\RR(\sigma)$ is contained in $\tau_i$ we obtain
$$P^\RR(\sigma) \; \subset \;  \tau_i\cap\tau_j \; \subset \; \tau_j.$$
By  Lemma \ref{gencone}, $\sigma$ is the convex hull of some rays
$\varrho_1, \ldots, \varrho_r$ of $\Delta$. For each of these rays we
have $P^\RR(\varrho_l) \subset \tau_j$. By the definition of
$\sigma_j$ all the rays $\varrho_l$ are contained in $\sigma_j$. This
implies $\sigma \subset {\sigma}_j$ and hence $\sigma = {\sigma}_i
\cap {\sigma}_j$. So $\sigma_i$ and $\sigma_j$ intersect in a face of
$\sigma_i$. That means $S$ is indeed a quasifan. Now we can apply
Remark \ref{quasifan} to conclude that
$G^\RR(\sigma_1),\dots,G^\RR(\sigma_r)$ in $\b{N}$ are the maximal
cones of the quotient fan $\b{\Delta}$  of $S$ by $L$. Moreover, $G$
defines a map from the fan $\Delta$ to the fan $\b{\Delta}$.
 
In the next step of the proof we show that $X_{\b{\Delta}}$ has a good
quotient by $H$. Since $L$ is contained in some $\sigma_i$ and $P^\RR$
maps $\sigma_i$ to the strictly convex cone $\tau_i$, we have $P(L) =
0$. It follows that $P$ defines an $L$-invariant map of systems of
cones from $S$ to $\t{\Delta}$. Consequently, there is a unique map of
fans $\b{P} : (\b{N}, \b{\Delta}) \to (\t{N}, \t{\Delta})$ with $P =
\b{P} \circ G$. Note that the associated toric morphism $\b{p} :
X_{\b{\Delta}} \to X \tq H$ is the toric quotient for $X_{\b{\Delta}}$
by $H$.

To check that $\b{p}$ is a good quotient, we use Proposition
\ref{goodcrit}. The first condition of \ref{goodcrit} is fulfilled
since by Remark \ref{maxcones}, $P^\RR(\sigma_i)=\tau_i$ for every $i$
and hence every maximal cone $\tau_i$ of $\t{\Delta}$ is the image
under $\b{P}^\RR$ of the maximal cone $G^\RR(\sigma_i)$ of
$\b{\Delta}$. For the verification of the second condition, let
$\b{\varrho}$ be a ray in $\b{\Delta}$ with $\b{P}^\RR(\b{\varrho})
\subset \tau_i$. Then by Lemma \ref{gencone} there is a ray
$\varrho\in\Delta$ with $G^\RR(\varrho)=\b{\varrho}$. Since
$P^\RR(\varrho)\subset \tau_i$, by definition $\varrho$ is contained
in $\sigma_i$ and hence $\b{\varrho}=G^\RR(\varrho)\subset
G^\RR(\sigma_i)$.

To complete the proof we have to verify the universal property of good
models  for $g$. So let ${X'}$ be a toric $H$-variety with a good
quotient $p' : X'  \to X' \tq H$ and let $f : X \to {X'}$ be a toric
$H$-morphism. Denote the fans associated to $X'$ and $X'\tq H$ by
$\Delta'$ and $\t{\Delta'}$ respectively, and let $F : (N,\Delta) \to
(N',\Delta')$ be the map of fans associated to $f$.

Now suppose for the moment that the linear map $F: N \to N'$ also
defines an $L$-invariant map from the system of cones $S$ to
$\Delta'$. Then, since $\b{\Delta}$ is the quotient fan of $S$ by $L$,
there is a unique map of fans $\b{F} : (\b{N}, \b{\Delta}) \to (N',
\Delta')$ with $F = \b{F} \circ G$. Clearly, the toric morphism $\b{f}
: X_{\b{\Delta}} \to X_{\Delta'}$ associated to $\b{F}$ provides us
with the required factorization of $f$ through $g$.

So it remains to show that for a given cone $\sigma_i \in S$ there is
a cone ${\sigma_i}'\in {\Delta'}$ with $F^\RR(\sigma_i) \subset
{\sigma_i}'$. (Since $L$ is contained in $\sigma_i$ and ${\sigma_i}'$
is strictly convex, this also implies that $F(L)=0$.) Consider the
following commutative diagrams of toric morphisms and of the
associated maps of fans:

\smallskip

$$\matrix{%
X & \bigtopmapright{f} &{X'} &\cr
\lmapdown{p} & & \rmapdown{p'} &\cr
X \tq H& \bigtopmapright{f'} & {X'}\tq H &, \cr}   \qquad \qquad \qquad
\matrix{%
(N,\Delta) & \bigtopmapright{F} & ({N'},{\Delta'}) &\cr
\lmapdown{P} & & \rmapdown{P'} &\cr
(\t{N},\t{\Delta}) & \bigtopmapright{F'} & (\t{N'},\t{\Delta'}) & .\cr}%
$$

\medskip

Let $\varrho$ be any ray of $\Delta$ which is contained in
$\sigma_i$. Since $P^\RR (\sigma_i) = \tau_i$, there is a maximal cone
$\tau_i'$ in $\t{\Delta'}$ containing  ${F'}^\RR
(P^\RR({\sigma}_i))={P'}^\RR(F^\RR({\sigma}_i))$. So in particular,
${P'}^\RR (F^\RR (\varrho)) \subset \tau_i'$. Suppose that $\sigma$ is
the minimal cone of $\Delta'$ containing $F^\RR(\varrho)$. Then
${P'}^\RR(\sigma)^\circ$ intersects $\tau_i'$ and therefore
${P'}^\RR(\sigma)$ is contained in $\tau_i'$.

Now $p'$ is a good quotient and therefore by Proposition
\ref{goodcrit}  there is a maximal cone ${\sigma_i}' \in \Delta'$ with
${P'}^\RR ({\sigma_i}') = \tau_i'$. Moreover, any cone of $\Delta'$
which is mapped into $\tau_i'$ by ${P'}^\RR$ is a face of
${\sigma_i}'$. So in particular, $\sigma \subset {\sigma_i}'$ and
hence $F^\RR(\varrho) \subset {\sigma_i}'$. Since $\sigma_i$ is
generated by the rays of $\Delta$ that it contains, we finally obtain
$F^\RR(\sigma_i) \subset {\sigma_i}'$. \endproof


\bibliography{}


\bigskip

{\sc

Fakult\"at f\"ur Mathematik und Informatik

Universit\"at Konstanz

Fach D197

D-78457 Konstanz

Germany}

\medskip

{\it E-mail address:}\enspace {\tt Annette.ACampo@uni-konstanz.de}

$\hphantom{\hbox{\it E-mail address:}}$\enspace {\tt Juergen.Hausen@uni-konstanz.de}

\end{document}